\documentclass{article}%
\usepackage{amsmath}
\usepackage{amsfonts}
\usepackage{amssymb}
\usepackage{graphicx}%
\providecommand{\U}[1]{\protect\rule{.1in}{.1in}}
\newtheorem{theorem}{Theorem}

\newtheorem{corollary}[theorem]{Corollary}

\newtheorem{example}[theorem]{Example}

\newtheorem{lemma}[theorem]{Lemma}

\newtheorem{proposition}[theorem]{Proposition}
\newtheorem{remark}[theorem]{Remark}

\newenvironment{proof}[1][Proof]{\noindent\textbf{#1.} }{\ \rule{0.5em}{0.5em}}
\begin{document}

\title{\textbf{Lifting diffeomorphisms} \\\textbf{to vector bundles}}
\author{\textsc{Jaime Mu\~{n}oz Masqu\'e }
\and \textsc{Mar\'{\i}a Eugenia Rosado Mar\'{\i}a}
\and \textsc{ Ignacio S\'anchez Rodr\'{\i}guez}}
\date{}
\maketitle

\begin{abstract}
\noindent Criteria for a diffeomorphism of a smooth manifold $M$ to be lifted
to a linear automorphism of a given real vector bundle $p\colon V\rightarrow
M$, are stated. Examples are included and the EEE metric and EEE complex
vector-bundle cases are also considered.

\end{abstract}

\smallskip



\smallskip

\noindent\emph{Acknowledgements:\/} M.E.R.M was partially supported by
Ministerio de Ciencia, Innovaci\'{o}n y Universidades (Spain) under Project no.\ PID2021--126124NB-I00.

\section{Introduction}

The goal of this paper is to study the problem of lifting a given
diffeomorphism $\varphi\colon M\rightarrow M$ of a connected, second countable
$C^{\infty}$ manifold $M$ without boundary of finite dimension to an
automorphism of a vector bundle $p\colon V\rightarrow M$ in the category of
$C^{\infty}$ vector bundles over $M$ (cf.\cite[\S 6]{KMS}). This question is
fundamental in the theory of transformation groups and topology.

The basic conditions given are as follows:

(1) the existence of a lift $\Phi$ only depends on the homotopy class of
$\varphi$,

(2) characteristic classes provide obstructions to the existence of a lift,

(3) $\varphi$ always lifts if the bundle is associated to a jet-automorphism
frame bundle.

The rest of the manuscript discusses whether a given lift can be enhanced with
further structure (metric, complex) or if it induces a lift to a frame bundle.

FIN PROVISIONAL

An automorphism of $V$ is a diffeomorphism $\Phi\colon V\to V$ permuting the
fibres of $p$, such that for every $x\in M$ the map $\Phi_{x}\colon V_{x}\to
V_{x^{\prime}}$ induced on the fibre $V_{x}=p^{-1}(x)$ is a linear map.

The group of diffeomorphisms of $M$ is denoted by $\mathrm{Diff\,}M$, and
below it is assumed to be endowed with the compact-open $C^{\infty}$-topology;
cf.\ \cite[\S 41.9]{KM}. For the consideration of other topologies on
$\mathrm{Diff\,}M$, we refer the reader to \cite{KM, KMR}.

Let $(\mathrm{AUT\,}V,\circ)$ be the group of all automorphisms of the vector
bundle $p\colon V\to M$. Each $\Phi\in\mathrm{AUT\,}V$ induces a
diffeomorphism $\varphi\colon M\to M$ on the base manifold by setting
$\varphi(x)=x^{\prime}$, where $x^{\prime}$ is given by $\Phi_{x}%
(V_{x})=V_{x^{\prime}}$, $\forall x\in M$, and the map $h_{V}\colon
\mathrm{AUT\,}V\to\mathrm{Diff\,}M$, $\Phi\mapsto\varphi$, is a group
homomorphism. Thus, the basic point is to decide under which hypotheses
$h_{V}$ is surjective. This problem arises in several situations; for example,
in studying Lie al\-ge\-bras of sections of a Lie algebra bundle and in
classifying vector bundles by their Lie algebras of linear endomorphism fields
(e.g., see \cite{L1}, \cite{L2}, \cite{L3}), and working with the group of all
automorphisms of a principal bundle and some of its subgroups (mainly the
gauge group) and its linear representations on vector bundles in gauge
theories (e.g., see \cite{B}, \cite{Bl}, \cite[III, \S \S 35, 36]{GS},
\cite{LR}, \cite{MM}, \cite{Mi}).

The difficulty of this lifting problem lies in that---unlike the gauge
group---the group of all automorphisms of $V$ cannot be written as the global
sections of a bundle over $M$ naturally attached to $V$, so that the problem
depends both on $V$ and $M$. These two groups (the full group of automorphisms
of $V$ and its gauge group) are defined in a precise way in the next section
\S \ref{prelim_def} below.

The infinitesimal version of the lifting problem (i.e., the lift of vector
fields) in the case of fibre bundles which are naturally attached to the first
order linear frame bundle---especially tensor bundles---has extensively been
dealt with in the literature (e.g., see \cite{LY}, \cite{Mo}, \cite{YL}) but
for vector fields and other operators such as connections, no topological
obstruction exists to lifting.

A summary of contents is as follows: In \S \ref{criteria} several criteria are
given so that the lifting problem to be possible (sufficient conditions for
lifting) and also some obstructions are obtained so that such a problem can be
carried out (necessary conditions to lifting). The obstructions to lifting are
expressed in terms of Stiefel-Whitney and Chern classes associated to the
vector bundle $p\colon V\rightarrow M$, depending on whether the vector bundle
is real or complex; see Proposition \ref{p0} and Corollary \ref{c0}.
Sufficient conditions to lifting are given in Proposition \ref{p1} and in
Proposition \ref{p2} the lifting problem is proved to depend only on the set
of homotopy classes of diffeomorphisms.

Section \ref{metric} is devoted to studying the particular case of a real
vector bundle $p\colon V\rightarrow M$ over a $C^{\infty}$-manifold, which is
assumed to be endowed with a positive definite symmetric bilinear form
$\langle\cdot,\cdot\rangle\colon V\oplus V=V\times_{M}V\rightarrow\mathbb{R}$.

Moreover, in \S \ref{complex} the lifting problem is studied in the category
of complex vector bundles over a connected, second countable $C^{\infty}$
manifold $M$ without boundary of finite dimension. Let $p\colon V\to M$ be a
complex vector bundle. In Proposition \ref{p4prime} the conditions are studied
so that a diffeomorphism that can be lifted to an $\mathbb{R}$-linear
automorphism of $V$, supports a $\mathbb{C}$-linear extension.

Finally, in \S \ref{examples}, several examples of relevant geometric interest
are studied in detail.

\section{Preliminaries and definitions}

\label{prelim_def}

EEE In this section we review quickly classical definitions on prinicpal and
vector bundles. In particular, we review the definition of the group of
vertical automorphisms of vector bundles, the $r$-th order frame bundle, the
gauge group and the fibre bundle associated with a principal $G$-bundle. EEE

The basic notations and definitions on principal and vector bundles have been
taken from \cite{Gr}, \cite{H}, \cite{KN}, \cite{K}, \cite{MS}, although below
we confine ourselves to exclusively deal with manifolds whose underlying
topology is second countable and, in addition, all manifolds considered are
assumed to be connected, without boundary, of finite dimension and of class
$C^{\infty}$, as well as the maps among them.

\subsection{The group $\mathrm{AUT}_{M}V$}

A \emph{morphism} between two vector bundles $p\colon V\rightarrow M$,
$p^{\prime}\colon V^{\prime}\rightarrow M^{\prime}$ is a smooth map
$\Phi\colon V\rightarrow V^{\prime}$ transforming each fibre of $p$ into a
fibre of $p^{\prime}$; i.e., $\Phi(V_{x})\subseteq V_{x^{\prime}}^{\prime}$,
in such a way that the restriction map $\Phi_{x}\colon V_{x}\rightarrow
V_{x^{\prime}}^{\prime}$ is $\mathbb{R}$-linear for all $x\in M$. The point
$x^{\prime}\in M^{\prime}$ satisfying the previous condition is unique and
hence, $\Phi$ induces a map $\varphi\colon M\rightarrow M^{\prime}$ by setting
$\varphi(x)=x^{\prime}$, which is differentiable by virtue of the
characteristic property of submersions and it is the only map such that,
$p^{\prime}\circ\Phi=\varphi\circ p$. The set of all morphisms from $V$ into
$V^{\prime}$ is denoted by $\mathrm{HOM}(V,V^{\prime})$. If $M=M^{\prime}$,
then the subset of all vertical (or fibred) morphisms (i.e., the morphisms
inducing the identity map on $M$) is denoted by $\mathrm{HOM}_{M}(V,V^{\prime
})$.

We remark that $\mathrm{HOM}_{M}(V,V^{\prime})$ is endowed with a natural
structure of a $C^{\infty}(M)$-module, for it is the module of global sections
of the vector bundle $\mathrm{Hom}_{M}(V,V^{\prime})\to M$, whose elements are
the linear maps $A\colon V_{x}\to V_{x}^{\prime}$, $x\in M$; i.e.,
$\mathrm{HOM}_{M}(V,V^{\prime})=\Gamma(M,\mathrm{Hom}_{M}(V,V^{\prime}))$.

Let $\mathrm{AUT\,}V$ (resp.\ $\mathrm{AUT}_{M}V$) be the subset of all
automorphisms (resp.\ vertical automorphisms) in $\mathrm{END\,}%
V=\mathrm{HOM\,}(V,V)$ (resp.\ $\mathrm{END}_{M}V=\mathrm{HOM}_{M}(V,V)$).
There is an exact sequence of groups
\[
1\rightarrow\mathrm{AUT}_{M}V\rightarrow\mathrm{AUT\,}V\overset{h_{V}%
}{\longrightarrow}\mathrm{Diff\,}M,\quad h_{V}(\Phi)=\varphi,
\]
which proves in particular that $\mathrm{AUT}_{M}V\subset\mathrm{AUT\,}V$ is a
normal subgroup. 

Let $\mathrm{Aut}_{M}V$ be the dense open subset in $\mathrm{End}_{M}V=
\mathrm{Hom}_{M}(V,V)$ of all automorphisms. The natural map
$\mathrm{Aut}_{M}V\rightarrow M$ induces a Lie-group fibre-bundle structure
over $\mathrm{Aut}_{M}V$ (cf.\ \cite{DL}) whose standard fibre is
$Gl(r,\mathbb{R})$, where $r=\operatorname*{rank}V$. The set of sections of
$\mathrm{Aut}_{M}V\rightarrow M$ is thus endowed with a group structure by
setting $(s\cdot s^{\prime})(x)=s(x)\cdot s^{\prime}(x)$, $\forall
s,s^{\prime}\in\Gamma(M,\mathrm{Aut}_{M}V)$, $\forall x\in M$, where, on the
right-hand side, the dot stands for the product in the Lie group
$(\mathrm{Aut}_{M}V)_{x}$. We thus have a natural group isomorphism
$(\mathrm{AUT}_{M}V,\circ)=(\Gamma(M,\mathrm{Aut}_{M}V),\cdot)$.

\subsection{The $r$-th frame bundle}

Let $\pi\colon F(M)\rightarrow M$ be the bundle of linear frames of $M$; its
points are ordered bases $(X_{1},\dotsc,X_{n})$ of the tangent space $T_{x}M$
at a point $x\in M$, $n=\dim M$. As is well known (e.g., see \cite[I,
\S 5]{KN}) $F(M)$ is a principal bundle with structure group the full linear
group $Gl(n,\mathbb{R})$. We define the $r$\emph{-th order frame bundle}
$\pi_{r}\colon F^{r}(M)\rightarrow M$ as the subset of all $r$-jets $j_{0}%
^{r}\xi\in J^{r}(\mathbb{R}^{n},M)$ such that $\xi_{\ast}\colon T_{0}%
\mathbb{R}^{n}\rightarrow T_{\xi(0)}M$ is a linear isomorphism, where $\pi
_{r}(j_{0}^{r}\xi)=\xi(0)$. Similarly to first order case, $F^{r}(M)$ is a
principal bundle with group $G^{r}(n)$, the $r$-th order jet group
(cf.\ \cite[IV.12]{KMS}); that is, the group of the $r$-jets $j_{0}^{r}f\in
J^{r}(\mathbb{R}^{n},\mathbb{R}^{n})$ such that $f(0)=0$ and $f_{\ast}\colon
T_{0}\mathbb{R}^{n}\rightarrow T_{0}\mathbb{R}^{n}$ is a linear isomorphism.

\subsection{Principal $G$-bundle, its associated bundles and the gauge group}

Let $G$ be a Lie group and let $\pi \colon P\rightarrow M$ be a principal
$G$-bundle (see\ \cite[Chapter 4, 2.2]{H}). If $G$ acts on the left on a 
manifold $F$, the \emph{fibre bundle associated with} $\pi\colon P\rightarrow
M$ and standard fibre $F$ (cf.\ \cite[I, \S 5]{KN}) is denoted by $\pi_{F}
\colon P\times^{G}F\rightarrow M$. If $F$ is a real or complex
finite-dimensional vector space and $G$ acts on $F$ through a linear
representation $\lambda\colon G\rightarrow Gl(F)$, then the associated bundle
is endowed with a natural structure of vector bundle.

An automorphism of $\pi \colon P\rightarrow M$ is a $G$-equivariant
diffeomorphism $\Phi \colon P\rightarrow P$. Every automorphisms $\Phi $ of
$P$, determines a unique diffeomorphism $\phi\colon M\rightarrow M$, such that
$\pi \circ \Phi =\phi \circ \pi $. If $\phi $ is the identity map on $M$, then
$\Phi $ is said to be a gauge transformation. The \emph{gauge group} of $P$ is
the group of all gauge transformations and is denoted by $\mathrm{Aut}_{M}P$ 
(cf.\ \cite[Chapter 7, 1]{H}).

Now, given a vector bundle $p\colon V\rightarrow M$, with $r=
\operatorname*{rank}V$, we obtain the \emph{linear frame bundle of} $V$ (cf.\ 
[13, III.10.11]) $\pi\colon FV \rightarrow M$ that is a principal $Gl(r,
\mathbb{R})$-bundle whose fibres are the ordered bases of $V_{x}$, $x\in M$.
The vector bundle $V$ itself is the bundle associated with $FV$ with the
natural action of $Gl(r,\mathbb{R})$ on $\mathbb{R}^{r}$. Also note that the
bundle $\mathrm{Aut}_{M}V\rightarrow M$ introduced above is naturally
isomorphic to the bundle associated with $FV$ with respect to the left adjoint
action of $Gl(r,\mathbb{R})$ into itself. Finally, it is straightforward
verify that the gauge group $\mathrm{Aut}_{M}FV$ is naturally isomorphic to
the group $\Gamma(M,\mathrm{Aut}_{M}V)$. 

\section{Basic criteria to lift}

\label{criteria}

NUEVO

\begin{proposition}
\label{Proposition1} If an automorphism $\varphi$ of $M$ can be lifted, then
$\varphi^{\ast}$ is the identity on the subring generated by the
characteristic classes.
\end{proposition}

\begin{proof}
If there exists an automorphism $\Phi\in\mathrm{AUT\,}V$ such that $h_{V}%
(\Phi)=\varphi$, then the vector bundles $V$ and $\varphi^{\ast}V$ are
$M$-isomorphic by means of the vertical isomorphism $V\rightarrow\varphi
^{\ast}V$, $v\mapsto(p(v),\Phi(v))$. From the properties of the
Stiefel-Whitney (resp.\ Chern) classes \emph{(cf.\ $SW_{1}$ and $C_{1}$ in
\cite[Chapter 17, \S 3.1 and \S 3.2]{H}) }we conclude that $\varphi$ keeps the
Stiefel-Whitney (resp.\ Chern) classes invariant, i.e., $H^{i}(\varphi
)(w_{i}(V))=w_{i}(V)$ (resp.\ $H^{2i}(\varphi)(c_{i}(V))=c_{i}(V)$), $1\leq
i\leq m$ (see $SW_{1}$ and $C_{1}$ in \cite[Chapter 17, \S 3.1 and \S 3.2]%
{H}). Therefore, $\varphi^{\ast}$ is the identity on the subring of $H^{\ast
}(M,\mathbb{Z}_{2})$ (resp. $H^{\ast}(M,\mathbb{R})$) generated by
Stiefel-Whitney (resp.\ Chern) classes.
\end{proof}

NUEVO

QUITAR

\begin{lemma}
\label{lemma1} \emph{(cf.\ $SW_{1}$ and $C_{1}$ in \cite[Chapter 17, \S 3.1
and \S 3.2]{H})} If $p_{i}\colon V_{i}\rightarrow M_{i}$ are real
(resp.\ complex) vector bundles for $i=1,2$, then it follows:

\emph{1)} If $M_{1}=M_{2}$ and $V_{1}$ and $V_{2}$ are $M$-isomorphic, then
$w_{j}(V_{1})=w_{j}(V_{2})$, (resp.\ $c_{j}(V_{1})=c_{j}(V_{2})$) for every
$0\leq j\leq m=\operatorname*{rank}V_{1}=\operatorname*{rank}V_{2}$.

\emph{2)} If $f\colon M_{1}\rightarrow M_{2}$ is a smooth map, then $f^{\ast
}(w_{j}V_{2})=w_{j}(f^{\ast}V_{2})$ (resp.\ $f^{\ast}(c_{j}V_{2}%
)=c_{j}(f^{\ast}V_{2})$).
\end{lemma}

\begin{proposition}
\label{p0} A necessary condition for a diffeomorphism $\varphi\colon
M\rightarrow M$ to be lifted to an automorphism of a real (resp.\ complex)
vector bundle $p\colon V\rightarrow M$ of rank $m$, is that $\varphi$ keeps
the Stiefel-Whitney (resp.\ Chern) classes invariant, i.e., $H^{i}%
(\varphi)(w_{i}(V))=w_{i}(V)$ (resp.\ $H^{2i}(\varphi)(c_{i}(V))=c_{i}(V)$),
$1\leq i\leq m$.
\end{proposition}

\begin{proof}
If there exists an automorphism $\Phi\in\mathrm{AUT\,}V$ such that $h_{V}%
(\Phi)=\varphi$, then the vector bundles $V$ and $\varphi^{\ast}V$ are
$M$-isomorphic by means of the vertical isomorphism $V\rightarrow\varphi
^{\ast}V$, $v\mapsto(p(v),\Phi(v))$, and we can conclude from the properties
of the Stiefel-Whitney (resp.\ Chern) classes stated in Lemma \ref{lemma1}.
\end{proof}

\begin{corollary}
\label{c0} If the Stiefel-Whitney (resp.\ Chern) classes of the vector bundle
$p\colon V\rightarrow M$ span $H^{\ast}(M,\mathbb{Z}_{2})$ (resp.\ $H^{\ast
}(M,\mathbb{R})$) and $\varphi\in\mathrm{Diff\,}M$ can be lifted to an
automorphism of $V$, then the homomorphism $\varphi^{\ast}\colon H^{\ast
}(M,\mathbb{Z}_{2})\rightarrow H^{\ast}(M,\mathbb{Z}_{2})$ (resp.\ $\varphi
^{\ast}\colon H^{\ast}(M,\mathbb{R})\rightarrow H^{\ast}(M,\mathbb{R})$)
induced in the cohomology ring, is the identity map.
\end{corollary}

QUITAR

Natural bundles are the bundles for which the lifting problem can be solved in
a functorial way. Each natural bundle over $n$-manifolds is equivalent to the
functor assigning the fibre bundle associated with $F^{r}(M)$ and standard
fibre $F$ to every $n$-dimensional manifold $M$, for certain integer $r$,
where $F$ is a fixed manifold on which $G^{r}(n)$ acts differentiably (see
\cite{ET}, \cite{PT}, \cite{T}). More precisely, a natural bundle is a
covariant functor from the category of $n$-dimensional manifolds and local
diffeomorphisms into the category of fibred manifolds and there exists a
bijective correspondence between the set of all $r$-th order natural bundles
on $n$-dimensional manifolds and the set of smooth left actions of the jet
group $G^{r}_{n}$ on smooth manifolds of finite dimension; see \cite[\S 14]%
{KMS}. Hence we have

\begin{proposition}
\label{p1} If $p\colon V\to M$ denotes the vector bundle associated with
$F^{r}(M)$ under a linear representation $\lambda\colon G^{r}(n)\to
Gl(m,\mathbb{R})$, then the homomorphism $h_{V}\colon\mathrm{AUT\,}%
V\to\mathrm{Diff\,}M$ is surjective.
\end{proposition}

Certainly, the criterion above does not suffice to solve the problem of
lifting diffeomorphisms to vector bundles; for example, in the category of
parallelizable manifolds each natural bundle is trivial. However, the lifting
problem can be posed for interesting examples of non-trivial bundles.

\begin{corollary}
\label{c1/2} If a diffeomorphism $\varphi\in\mathrm{Diff\,}M$ exists that
cannot be lifted to the vector bundle $p\colon V\to M$, then $V$ is not
associated under a linear representation to $F^{r}(M)$ for any integer
$r\in\mathbb{N}$. In particular, $V$ is not isomorphic to a tensorial bundle
$(\otimes^{p}T^{\ast}M)\otimes(\otimes^{q}TM)$, $p,q\in\mathbb{N}$.
\end{corollary}

\begin{proposition}
\label{p2} If a diffeomorphism $\varphi\colon M\to M$ is homotopic to the
identity map, then there exists $\Phi\in\mathrm{AUT\,}V$ such that $h_{V}%
(\Phi)=\varphi$.
\end{proposition}

\begin{proof}
We know (e.g., see \cite[IV, 4.11 Corollaire]{G}) that if two differentiable
maps are $C^{0}$ homotopic, then they are also $C^{\infty}$ homotopic; i.e.,
there exists a differentiable mapping $\phi\colon\lbrack0,1]\times M\to M$
such that $\phi_{0}=\mathrm{id}_{M}$, $\phi_{1}=\varphi$. Now we use the
parallel displacement along curves in vector bundles with respect to a
connection $\nabla$ on $V$. If $\sigma\colon\lbrack0,1]\to M$ is a curve and
$v_{0}\in V_{x_{0}}$ is a given vector over $x_{0}=\sigma(0)$, then a unique
smooth section $\chi\colon\lbrack0,1]\to V$ exists such that, 1) $\chi$ is
defined along $\sigma$, i.e., $\chi(t)\in V_{\sigma(t)}$, $\forall t\in
\lbrack0,1]$, 2) $\chi(0)=v_{0}$, and 3) $\nabla_{T}\chi=0$, where $T$ is the
tangent vector to $\sigma$.

Hence we have a linear isomorphism $\tau_{\sigma}\colon V_{\sigma(0)}$ $\to
V_{\sigma(1)}$ (the parallel displacement) by setting $\tau_{\sigma}%
(v_{0})=\chi(1)$ for every $v_{0}\in V_{\sigma}(0)$, where $\chi$ is as above
and we define $\Phi$ as follows: If $v\in V_{x}$ then $\Phi(v)=\tau_{\sigma
}(v)$, where $\sigma$ is the curve given by $\sigma(t)=\phi_{t}(x)$, $0\leq
t\leq1$. It is easy to check that $\Phi$ is an automorphism of $V$, and
certainly $h_{V}(\Phi)=\varphi$, as $\sigma(1)=\phi_{1}(x)=\varphi(x)$.
\end{proof}

\begin{corollary}
\label{c1} If $\varphi$, $\psi\colon M\to M$ are two $C^{0}$ homotopic
diffeomorphisms, then $\varphi$ can be lifted to an automorphism of $V$ if and
only if $\psi$ can. Hence if any diffeomorphism of $M$ is homotopic to a
diffeomorphism of $M$ that can be lifted to an automorphism of $V$, then
$h_{V}\colon\mathrm{AUT\,}V\to\mathrm{Diff\,}M$ is surjective.
\end{corollary}

\begin{proof}
If $\phi_{t}\colon M\to M$ is a $C^{0}$ homotopy from $\varphi$ to $\psi$,
then $\varphi^{-1}\circ\phi_{t}\colon M\to M$ is a homotopy from
$\text{Id}_{M}$ to $\varphi^{-1}\circ\psi$. Therefore there exists an
automorphism $\Upsilon\in\mathrm{AUT\,}V$ such that $h_{V}(\Upsilon
)=\varphi^{-1}\circ\psi$. Since $\varphi$ lifts to an automorphism $\Phi
\in\mathrm{AUT}\,V$, then $h_{V}(\Phi)=\varphi$ and $h_{V}(\Phi\circ
\Upsilon)=\varphi\circ\varphi^{-1}\circ\psi=\psi$; hence $\psi\in
\mathrm{Diff\,}M$ lifts to $\Psi=\Phi\circ\Upsilon\in\mathrm{AUT\,}V$.
\end{proof}

Let $\mathrm{Diff}_{1}M$ be the normal subgroup in $\mathrm{Diff\,}M$ of
diffeomorphisms homotopic to the identity map of $M$. According to Proposition
\ref{p2}, if the equivalence class $c(\varphi)\in\mathrm{Diff\,}%
M/\mathrm{Diff}_{1}M$ of a diffeomorphism $\varphi\in\mathrm{Diff\,}M$
vanishes, then $\varphi$ can be lifted to an automorphism of the vector bundle
$p\colon V\to M$. Therefore, the obstruction to lifting lies on this class.

\begin{remark}
\label{remark1} Let $\mathrm{Diff}_{0}M$ denote the normal subgroup in
$\mathrm{Diff\,}M$ of diffeomorphisms isotopic to the identity map of $M$. We
have $\mathrm{Diff}_{0}M\subseteq\mathrm{Diff}_{1}M$, and the quotient group
$\mathrm{Diff}_{1}M/\mathrm{Diff}_{0}M$ is a discrete group, as $\mathrm{Diff}%
_{0}M$ is a path component of $\mathrm{Diff\,}M$; hence also a connected
component since $\mathrm{Diff\,}M$ is locally path-cconnected, and
$\mathrm{Diff}_{1}M$ is a union of components of $\mathrm{Diff\,}M$.
Nevertheless, in general $\mathrm{Diff}_{0}M$ and $\mathrm{Diff}_{1}M$ are
distinct; e.g., see \cite{Ruber}.
\end{remark}

\begin{remark}
\label{remark2} As statd above (Proposition \ref{p2}), the lifting problem
depends only on the set of homotopy classes of diffeomorphisms, although this
set may be non-trivial (and, in fact, quite large; for example, for the
$n$-dimensional torus $(S^{1})^{n}$ this set can be identified to the matrices
$A\in Gl(n,\mathbb{Z})$ such that $\det A=\pm1$).

Certainly, the lifting problem can immediately be solved for some classes of
vector bundles; namely, for natural vector bundles (\cite{ET}, \cite{PT},
\cite{T}) as shown in Proposition \ref{p1} above, but the difficulty
essentially subsists since no criterion is known in order to decide whether a
vector bundle, given in advance, is natural or not.
\end{remark}

\begin{example}
\label{Ex1} Given a vector bundle $p\colon V\to M$, for a large enough $N$,
there exists an $M$-monomorphism of vector bundles $\iota\colon V\to
M\times\mathbb{R}^{N}$. Let $p^{\prime}\colon V^{\prime}\to M$ be the vector
bundle orthogonal to $V$ in $\mathbb{R}^{N}$, i.e., $\mathbb{R}^{N}%
=\iota(V_{x})\oplus V_{x}^{\prime}$, $\forall x\in M$, and let $\varpi\colon
M\times\mathbb{R}^{N}\to\iota(V)$ be the projection onto the first summand.
Every $\varphi\in\mathrm{Diff\,}M$ admits an obvious lift $(\varphi
,\mathrm{id}_{\mathbb{R}^{N}})$ to $M\times\mathbb{R}^{N}$. If the composite
mapping $\varpi\circ(\varphi,\mathrm{id}_{\mathbb{R}^{N}})\circ\iota\colon
V\to\iota(V)$ is injective, then it determines a lift $\Phi$ of $\varphi$ to
$V$ such that, $\varpi\circ(\varphi,\mathrm{id}_{\mathbb{R}^{N}}) \circ
\iota=\iota\circ\Phi$.
\end{example}

\section{Reduction to metric bundles}

\label{metric}

\begin{proposition}
\label{p3} Let $p\colon V\to M$ be a real vector bundle over a $C^{\infty}%
$-manifold, which is assumed to be endowed with a positive definite symmetric
bilinear form $\langle\cdot,\cdot\rangle\colon V\oplus V =V\times_{M}%
V\to\mathbb{R}$. A diffeomorphism $\varphi\colon M\to M$ can be lifted to an
automorphism $\Phi\colon V\to V$ if and only if $\varphi$ can be lifted to an
isometry of the metric vector bundle $(V,\langle\cdot,\cdot\rangle)$.
\end{proposition}

\begin{proof}
If $\varphi$ can be lifted to $\Phi\colon V\to V$, then from the polar
decomposition there exist unique linear maps $\Psi_{x}\colon V_{x}\to
V_{\varphi(x)}$, $\Xi_{x}\colon V_{\varphi(x)}\to V_{\varphi(x)}$ such that
for every $x\in M$,

\begin{enumerate}
\item[(i)] $\Psi_{x}\colon(V_{x},\langle\cdot,\cdot\rangle_{x}) \to
(V_{\varphi(x)},\langle\cdot,\cdot\rangle_{\varphi(x)})$ is an isometry,

\item[(ii)] $\Xi_{x}\colon(V_{\varphi(x)},\langle\cdot,\cdot$ $\rangle
_{\varphi(x)}) \to(V_{\varphi(x)},\langle\cdot,\cdot$ $\rangle_{\varphi(x)})$
is a positive definite symmetric automorphism, and

\item[(iii)] $\Phi_{x}=\Xi_{x}\circ\Psi_{x}$.
\end{enumerate}

Note that the polar decomposition of $\Phi$ depends smoothly on $\Phi$. Hence
we obtain a vector bundle linear isometry $\Psi\colon(V,\langle\cdot
,\cdot\rangle)\to$ $(V,\langle\cdot,\cdot\rangle)$ that lifts $\varphi$. The
converse is trivial.
\end{proof}

\begin{corollary}
\label{c2} Let $g$ be an arbitrary Riemannian metric on $M$. Every
diffeomorphism $\varphi\in\mathrm{Diff\,}M$ can be lifted to an isometry of
$\left(  V,\langle\cdot,\cdot\rangle\right)  =(TM,g)$.
\end{corollary}

\begin{proof}
Let $\pi\colon F(M)\to M$ be the bundle of linear frames of $M$. As the
tangent bundle $TM$ is an associated bundle to $F(M)$ (e.g., see \cite[I,
Example 5.3]{KN}), from Proposition \ref{p1} the diffeomorphism $\varphi$ can
be lifted to $TM$, and the result thus follows from Proposition \ref{p3}.
\end{proof}

\begin{proposition}
\label{p4} Assume $p\colon V\to M$ and $\langle\cdot,\cdot\rangle\colon
V\oplus V=V\times_{M}V\to\mathbb{R}$ are as in \emph{Proposition \ref{p3}},
with $r=\operatorname*{rank}V$. Let $\pi\colon OF(V)\to M$ be the principal
$O(r)$-bundle of orthonormal frames of $V$. A diffeomorphism $\varphi\colon
M\to M$ can be lifted to an isometry of the metric vector bundle
$(V,\langle\cdot,\cdot\rangle)$ if and only if $\varphi$ can be lifted to an
automorphism of the principal bundle $\pi\colon OF(V)\to M$. If $p\colon V\to
M$ is an oriented vector bundle, then a similar result holds for the principal
$SO(r)$-bundle of positively orthonormal frames $\pi\colon O^{+}F(V)\to M$ of
$V$.
\end{proposition}

\begin{proof}
If $\varphi$ lifts to an isometry $\Phi\colon V\to V$, then we can define a
principal bundle automorphism $\Psi\colon OF(V)\to OF(V)$ by simply setting
\[
\Psi(v_{1},\dotsc,v_{r})=(\Phi(v_{1}),\dotsc,\Phi(v_{r})), \quad\forall
(v_{1},\dotsc,v_{r})\in OF(V),
\]
as $\Phi$ transforms an orthonormal frame into another. The converse follows
taking account of the fact that $V$ is isomorphic to the bundle associated
with $OF(V)$ and $\mathbb{R}^{r}$ with respect to the standard action of the
orthogonal group on vectors; i.e., $V=OF(V)\times^{O(r)}\mathbb{R}^{r}$.
\end{proof}

As is well-known (e.g., see \cite{Gr}), $C^{\infty}$ real line bundles over
$M$ are classified by the group $H^{1}(M;\mathbb{Z}_{2})$ and there is a
natural and one-to-one correspondence between real line bundles and $2$-fold
coverings on a manifold. Hence by applying Propositions \ref{p3}, \ref{p4},
and the principle of monodromy (e.g., see \cite[IX, 3.5 Proposition]{G}) we obtain

\begin{corollary}
\label{c3} Let $p\colon\tilde{M}_{L}\to M$ be the $2$-fold covering associated
to the real line bundle $L\to M$. A diffeomorphism $\varphi\colon M\to M$ can
be lifted to an automorphism of $L$ if and only if $\varphi_{\ast}\colon
\pi_{1}(M)\to\pi_{1}(M)$ maps the subgroup $p_{\ast}\pi_{1}(\tilde{M}_{L})$
onto itself. Therefore, for real line bundles the lifting problem is reduced
to a topological question. In particular, every diffeomorphism of $M$ can be
lifted to any real line bundle if $\pi_{1}(M)$ has a unique subgroup of index
$2$.
\end{corollary}

\begin{example}
\label{Ex2} Let $p\colon V\to M$ be a real vector bundle, let $\varphi$ be a
diffeomorphism of $M$, and let $U_{1}$, $U_{2}$ be two open subsets of $M$
such that, i) $\varphi(U_{i})=U_{i}$, $i=1,2$, and ii) $p\colon p^{-1}%
(U_{i})\to U_{i}$, $i=1,2$, is trivial. Let $s_{i}\colon U_{i}\to OF(V)$,
$i=1,2$, be two smooth sections. According to the local structure of the gauge
group (e.g., see \cite[Chapter 7, 1.7 Proposition]{H}), every lift of
$\varphi|_{U_{i}}$ to $\pi^{-1}(U_{i})$ can be written as $\Phi_{i}%
(s_{i}(x)\cdot A_{i})=s_{i}(\varphi(x))\cdot\alpha_{i}(x)A_{i}$, $\forall x\in
U_{i}$, $\forall A_{i}\in O(r)$, and for certain $\alpha_{i}\in C^{\infty
}(U_{i},O(r))$, $i=1,2$. Moreover, there exists $\alpha\in C^{\infty}%
(U_{1}\cap U_{2},O(r))$ such that $s_{2}(x)=s_{1}(x)\cdot\alpha(x)$, $\forall
x\in U_{1}\cap U_{2}$, and $\varphi$ can be lifted to $U_{1}\cup U_{2}$ if and
only if, $\Phi_{1}|_{U_{1}\cap U_{2}}=\Phi_{2}|_{U_{1}\cap U_{2}}$ or
equivalently, $\alpha_{2}(x)=\alpha(\varphi(x))^{-1}\alpha_{1}(x)\alpha(x)$,
$\forall x\in U_{1}\cap U_{2}$.
\end{example}

\section{The complex case}

\label{complex}

\begin{proposition}
\label{p4prime} Let $p\colon V\to M$ be a complex vector bundle of complex
rank $r$. If a diffeomorphism $\varphi\colon M\to M$ can be lifted to a
$\mathbb{R}$-linear automorphism $\Phi\colon V\to V$, then there exists a
fibre bundle $q\colon E\to M$ with standard fibre $Gl(r,\mathbb{C})$ whose
global sections can be identified to the $\mathbb{C}$-linear automorphisms
$\Phi^{\prime}\colon V\to V$ such that $h_{V}(\Phi^{\prime})=\varphi$.
\end{proposition}

\begin{proof}
If $\Phi^{\prime}$ is another $\mathbb{R}$-linear automorphism of $V$ such
that $h_{V}(\Phi^{\prime})=\varphi$, then $\Psi=\Phi^{\prime}\circ\Phi^{-1}$
is a vertical automorphism, i.e., $\Psi\in\mathrm{AUT}_{M}V$, and
$\Phi^{\prime}$ is $\mathbb{C}$-linear if and only if $\Phi^{\prime}$ commutes
with the endomorphism $J\colon V\to V$ defining the complex structure on $V$;
i.e., $J\circ\Psi\circ\Phi=\Psi\circ\Phi\circ J$, or equivalently,
\[
\Psi_{x}^{-1}\circ J_{x}\circ\Psi_{x} =\Phi_{\varphi^{-1}(x)}\circ
J_{\varphi^{-1}(x)} \circ(\Phi_{\varphi^{-1}(x)})^{-1},\quad\forall x\in M.
\]
Note that $(\Phi^{-1})_{x}=(\Phi_{\varphi^{-1}(x)})^{-1}$. We define $K\colon
V\to V$ by setting,
\[
K_{x}=\Phi_{\varphi^{-1}(x)}\circ J_{\varphi^{-1}(x)} \circ(\Phi_{\varphi
^{-1}(x)})^{-1},
\]
and for every $x\in M$, let $E_{x}$ be the set of $\mathbb{R}$-linear
automorphisms $\psi\colon V_{x}\to V_{x}$ such that $\psi^{-1}\circ J_{x}%
\circ\psi=K_{x}$. It suffices to prove that $E=\coprod_{x\in M}E_{x}$ is a
fibre bundle with respect to the natural projection $q\colon E\to M$, whose
standard fibre is $Gl(r,\mathbb{C})$. Note that $E_{x}$ is non-empty for every
$x\in M$, because $\Psi_{x}\in E_{x}$.

Once a point $x_{0}\in M$ has been fixed, let $U$ be an open neighbourhood of
$x_{0}$ in $M$ such that $V$ is trivial over $U$ as a complex vector bundle
and let $v_{j}\colon U\to V$, $1\leq j\leq r$ be a basis for the module of sections.

If $\psi,\psi^{\prime}\in E_{x}$, then $\lambda=\psi^{\prime}\circ\psi^{-1}$
commutes with $J_{x}$. Hence the matrix of $\lambda$ in the basis $\left(
v_{i}(x),J(v_{i}(x))\right)  _{i=1}^{r}$ takes the form
\[
\Lambda=\left(
\begin{array}
[c]{rr}%
A & -B\\
B & A
\end{array}
\right)  ,\qquad A+Bi\in Gl(r,\mathbb{C}).
\]
Hence $\psi^{\prime}=\lambda\circ\psi$, thus finishing the proof.
\end{proof}

\begin{corollary}
\label{c4} With the same assumptions as in \emph{Proposition \ref{p4prime}},
if we further assume $H^{k}(M;\pi_{k-1}(U(r)))=0$ for $1\leq k\leq n$, then
there exists a $\mathbb{C}$-linear automorphism of $V$ lifting $\varphi
\in\mathrm{Diff\,}M$.
\end{corollary}

\begin{proof}
This is an easy consequence of Proposition \ref{p4prime} and the obstruction
theory (the argument is the same as in \cite[IV,~section 44, pp.~351]{GS}),
taking the diffeomorphism $Gl(r,\mathbb{C})\cong U(r)\times\mathbb{R}^{r^{2}}$
into account.
\end{proof}

Given a complex vector bundle $p\colon V\to M$, we say that the $\mathbb{C}%
$-linear automorphism lifting problem can be solved for $V$ if every
$\varphi\in\mathrm{Diff\,}M$ admitting a $\mathbb{R}$-linear automorphism
$\Phi\in\mathrm{AUT\,}V$ such that $h_{V}(\Phi)=\varphi$, also admits a
$\mathbb{C} $-linear automorphism lifting it to $V$.

\begin{corollary}
\label{c5} We have

\begin{enumerate}
\item[\emph{(i)}] On a connected non-compact surface $M$, the $\mathbb{C}%
$-linear automorphism lifting problem can always be solved.

\item[\emph{(ii)}] On a $3$-dimensional manifold $M$ the $\mathbb{C}$-linear
automorphism lifting problem can be solved if $H^{2}(M;\mathbb{Z})=0$.
\end{enumerate}
\end{corollary}

\begin{proof}
Item (i) follows from $\pi_{1}(U(r))=\mathbb{Z}$ (see \cite[Chapter 8,
\S 12.3]{H}) and Corollary \ref{c4}, taking account of the fact that
$H^{2}(M;\mathbb{Z})=0$ by virtue of the hypothesis that $M$ is a connected
non-compact surface. Similarly, (ii) follows from $\pi_{2}(U(r))=0$ (see
\cite[Chapter 8, \S 12.4]{H}) and Corollary \ref{c4}.
\end{proof}

\section{Some examples}

\label{examples}

\begin{example}
\label{Ex3} As the Stiefel-Whitney (resp.\ Chern) classes of the tautological
vector bundle $\gamma^{k}(\mathbb{R}^{n})\to$ $G_{k}(\mathbb{R}^{n})$
(resp.\ $\gamma^{k}(\mathbb{C}^{n})\to$ $G_{k}(\mathbb{C}^{n})$) span
$H^{\ast}(G_{k}(\mathbb{R}^{n}),\mathbb{Z}_{2})$ (resp.\ $H^{\ast}%
(G_{k}(\mathbb{C}^{n}),\mathbb{R})$), if a diffeomorphism $\varphi$ of
$G_{k}(\mathbb{R}^{n})$ (resp.\ $G_{k}(\mathbb{C}^{n})$) can be lifted to an
automorphism of $\gamma^{k}(\mathbb{R}^{n})$ (resp.\ $\gamma^{k}%
(\mathbb{C}^{n})$), then the homomorphism induced in the cohomology ring,
$\varphi^{\ast}\colon H^{\ast}(G_{k}(\mathbb{R}^{n}),\mathbb{Z}_{2}) \to
H^{\ast}(G_{k}(\mathbb{R}^{n}),\mathbb{Z}_{2})$ (resp.\ $\varphi^{\ast}\colon
H^{\ast}(G_{k}(\mathbb{C}^{n}),\mathbb{R}) \to H^{\ast}(G_{k}(\mathbb{C}%
^{n}),\mathbb{R})$), is the identity map. In the complex case, this means that
the homotopy class of $\varphi$ belongs to the group $\mathcal{E}^{\ast}%
(G_{k}(\mathbb{C}^{n}))$; see \cite{ArkMaru}, \cite{Baues}, \cite{Pav}. (For
maps between distinct complex Grassmannians of the same dimension, the
situation seems to be rather different, see \cite[Theorem 1.3]{SS}.) From
Corollary \ref{c1/2} it follows that the tautological bundle $\gamma
^{k}(\mathbb{C}^{n})$ is \emph{not} a tensorial bundle.
\end{example}

\begin{example}
\label{Ex4} If a diffeomorphism $\varphi$ of the Grassmannian $G_{k}%
(\mathbb{R}^{n})$ of $k$-planes into itself is such that, $S\cap\varphi(S)=\{
0\} $ for every $S\in G_{k}(\mathbb{R}^{n})$ (which implies $k\leq\frac{1}%
{2}n$), then $\varphi$ can be lifted to an automorphism $\Phi$ of the
tautological vector bundle $\gamma^{k}(\mathbb{R}^{n})\to$$G_{k}%
(\mathbb{R}^{n})$ (cf.\ \cite{MS}) by setting $\Phi(S,v)=(\varphi
(S),v_{\varphi(S)}^{\bot})$, where $v_{S}^{\bot}$ denotes the orthogonal
projection of the vector $v$ onto the subspace $S^{\bot}$. In particular, this
applies to the canonical involution $\varphi\colon G_{k}(\mathbb{R}^{2k})\to
G_{k}(\mathbb{R}^{2k})$, $\varphi(S)=S^{\bot}$.
\end{example}

\begin{example}
\label{Ex4a} As is well known, every symmetric idempotent linear mapping
$P\colon\mathbb{R}^{n}\to\mathbb{R}^{n}$ of trace equal to $k$ is the
orthogonal projection onto the subspace $S_{P}=\{ v\in\mathbb{R}^{n}:P(v)=v\}
$. Hence, every $A\in O(n)$ induces a diffeomorphism $\varphi_{A}\colon
G_{k}(\mathbb{R}^{n})\to G_{k}(\mathbb{R}^{n})$, $\varphi_{A}(S_{P}%
)=S_{APA^{t}}$, which admits a natural lift to the tautological bundle
$\Phi_{A}\colon\gamma^{k}(\mathbb{R}^{n})\to\gamma^{k}(\mathbb{R}^{n})$ by
setting $\Phi_{A}(S_{P},v)=(\varphi_{A}(S_{P}),A(v))$ for all $(S_{P}%
,v)\in\gamma^{k}(\mathbb{R}^{n})$.
\end{example}

\begin{example}
\label{Ex4.5} As a consequence of Corollary \ref{c3}, we deduce that every
diffeomorphism of the real projective space can be lifted to an automorphism
of the tautological line bundle $\gamma^{1}(\mathbb{R}^{n+1})\to
\mathbb{R}P^{n}=G_{1}(\mathbb{R}^{n+1})$ as the $2$-fold covering associated
to $\gamma^{1}(\mathbb{R}^{n+1})$ is $S^{n}\to\mathbb{R}P^{n}$. Recall that
the fibre bundle $S^{n}\to\mathbb{R}P^{n}$ is universal for dimensions $\leq
n-1$ (\cite[Chapter 4, 11.3]{H}).
\end{example}

\begin{example}
\label{Ex5} Every diffeomorphism $\varphi$ of the complex projective space
$\mathbb{C}P^{n}$ can be lifted to a $\mathbb{R}$-linear automorphism of the
complex tautological line bundle,
\[
\gamma^{1}(\mathbb{C}^{n+1})\to\mathbb{C}P^{n} =G_{1}(\mathbb{C}^{n+1}).
\]
In fact, we know (see \cite[Chapter 8, \S 4, 11 Theorem]{Sp}) that if $X$ is a
$CW$ complex such that, 1) $\dim X\leq2n+1$, 2) $H^{2n+1}(X;\mathbb{Z})=0$,
then the homotopy classes $[X,\mathbb{C}P^{n}]$ are classified by
$H^{2}(X;\mathbb{Z})$ under the map $[f]\mapsto f^{\ast}(\iota)$, where
$\iota$ is a generator of $H^{2}(\mathbb{C}P^{n};\mathbb{Z})$. Letting
$X=\mathbb{C}P^{n}$, we conclude that there are two homotopy classes of
diffeomorphisms from $\mathbb{C}P^{n}$ onto itself: That of the identity map
and that of $\varphi\colon\mathbb{C}P^{n} \to\mathbb{C}P^{n}$, $\varphi
(x)=\bar{x}$, where $\bar{x}$ denotes the conjugate point; i.e., if the
homogeneous coordinates of $x$ are $[z_{0},\dotsc,z_{n}]$, then the
homogeneous coordinates of $\bar{x}$ are $[\bar{z}_{0},\dotsc,\bar{z}_{n}]$.
According to Corollary \ref{c1}, we only need to lift $\varphi$. This is done
as follows. A point in $\gamma^{1}(\mathbb{C}^{n+1})$ is a pair $(x=[z_{0}%
,\dotsc,z_{n}],w)$ where $x\in\mathbb{C}P^{n}$ and $w\in\mathbb{C}^{n+1}$ is a
vector such that $w=\lambda\cdot(z_{0},\dotsc,z_{n})$, with $\lambda
\in\mathbb{C}$. We set $\Phi(x,w)=(\bar{x},\bar{\lambda} \cdot(\bar{z}%
_{0},\dotsc,\bar{z}_{n}))$. The definition is easily seen to be independent of
the homogeneous coordinates chosen for the point $x$, thus providing the
desired lifting. Similarly to Example \ref{Ex4.5}, if $\left\langle
\cdot,\cdot\right\rangle $ denotes the standard Hermitian product on
$\gamma^{1}(\mathbb{C}^{n+1})$, then the Hopf fibration $S^{2n+1}\to
\mathbb{C}P^{n}$ coincides with the principal $SO(2)$-bundle of positively
oriented orthonormal frames of $\gamma^{1}(\mathbb{C}^{n+1})$
(cf.\ Proposition \ref{p4}).

The case of the quaternionic tautological bundle $\gamma^{1}(\mathbb{H}%
^{n+1})\to\mathbb{H}P^{n}$ and its associated Hopf fibration $S^{4n+3}%
\to\mathbb{H}P^{n}$ can be dealt with similarly.
\end{example}

\begin{example}
\label{Ex6} According to Hopf classification theorem (e.g., see \cite[Chapter
8, \S 1, 16 Corollary]{Sp}) and Corollary \ref{c1} above, every $\varphi
\in\mathrm{Diff\,}S^{n}$ of degree $+1$ can be lifted to any vector bundle
over $S^{n}$. Next, every $\varphi\in\mathrm{Diff\,}S^{2}$ of degree $-1$ on
the $2$-sphere is proved to admit a lift to an arbitrary vector bundle. In
this case, $\varphi$ is homotopic to the symmetry $\sigma\colon S^{2}\to
S^{2}$, $\sigma(x_{1},x_{2},x_{3})=(x_{1},-x_{2},x_{3})$, which, in
homogeneous coordinates (taking the usual isomorphism $S^{2}\cong%
\mathbb{C}P^{1}$ into account) transforms into $\sigma([z_{0},z_{1}])=[\bar
{z}_{0},\bar{z}_{1}]$, $\forall\lbrack z_{0},z_{1}]\in\mathbb{C}P^{1}$. We can
confine ourselves to consider vector bundles of rank $2$ on $S^{2}$. As is
known (e.g., see \cite{Moore}) each of such bundles is isomorphic to
$(f_{n})^{\ast}\gamma^{1}(\mathbb{C}^{2})$, $n\in\mathbb{N}$, where
$f_{n}\colon\mathbb{C}P^{1}\to\mathbb{C}P^{1}$ is the mapping $f_{n}%
([z,1])=[z^{n},1]$, $f_{n}([1,0])=[1,0]$. By proceeding as in Example
\ref{Ex5}, a real lift $\Phi_{n}\colon(f_{n})^{\ast}\gamma^{1}(\mathbb{C}^{2})
\to(f_{n})^{\ast}\gamma^{1}(\mathbb{C}^{2})$ of $\sigma$ is obtained by
setting $\Phi_{n} \left(  [z_{0},z_{1}], \lambda\left(  (z_{0})^{n}%
,(z_{1})^{n} \right)  \right)  =\left(  [\bar{z}_{0},\bar{z}_{1}],
\bar{\lambda} \left(  (\bar{z}_{0})^{n},(\bar{z}_{1})^{n} \right)  \right)  $.
\end{example}

\begin{example}
According to Propositions \ref{p1}, \ref{p3}, and \ref{p4}, every $\varphi
\in\mathrm{Diff\,}S^{n}$ can be lifted to the principal $O(n)$-bundle
$O(n+1)\to S^{n}$, as this bundle can be identified to the bundle of
orthonormal linear frames on $S^{n}$ with respect to the Euclidean metric.
\end{example}

\begin{example}
For every sequence $b=(b_{1},\dotsc,b_{n})\in\{ 0,1\} ^{n}$, $b\neq
(0,\dotsc,0)$, let $S_{b}\subset\mathbb{Z}^{n}$ be the subgroup of the
elements $(x_{1},\dotsc,x_{n})\in\mathbb{Z}^{n}$ for which the integer
$b_{1}x_{1}+\dotsc+b_{n}x_{n}$ is even and let $L_{b}\to\mathbb{R}%
^{n}/\mathbb{Z}^{n}$ be the line bundle corresponding to the $2$-fold covering
$\mathbb{R}^{n}/S_{b}\to\mathbb{R}^{n}/\mathbb{Z}^{n}$. With the natural
identification $\pi_{1} \left(  \mathbb{R}^{n}/\mathbb{Z}^{n} \right)
\cong\pi_{1}(\mathbb{R}/\mathbb{Z})^{n}\cong\mathbb{Z}^{n}$, a diffeomorphism
$\varphi\colon\mathbb{R}^{n}/\mathbb{Z}^{n}\to\mathbb{R}^{n}/\mathbb{Z}^{n}$
can be lifted to an automorphism of $L_{b}$ if and only if $\pi_{1}%
(\varphi)(S_{b})\subseteq S_{b}$.
\end{example}

\begin{example}
According to \cite{DG}, the diffeomorphisms $a,r,s\colon S^{1}\times S^{2}\to
S^{1}\times S^{2}$, $a(\theta,z)=(\theta,-z)$; $r(\theta,x^{1}+ix^{2}%
,x^{3})=(\theta,\exp(i\theta)(x^{1}+ix^{2}),x^{3})$, $\sum_{i=1}^{3}%
(x^{i})^{2}=1$; $s(\theta,z)=(-\theta,z)$, span the diffeotopy group of
$S^{1}\times S^{2}$. Below, we use the usual identifications: $S^{1}%
=\mathbb{R}P^{1}$, $S^{2}=\mathbb{C}P^{1}$. With the same notations as in
Example \ref{Ex6}, the diffeomorphisms $a$ and $s$ can be lifted to the
product bundle $\gamma^{1}(\mathbb{R})\times(f_{n})^{\ast}\gamma
^{1}(\mathbb{C}^{2}) \to S^{1}\times S^{2}$. In fact, if $\sigma^{1}\colon
S^{1}\to S^{1}$ is the conjugation map, i.e., $\sigma^{1}(z)=\bar{z}$ and
$\sigma^{2}\colon S^{2}\to S^{2}$ is the symmetry $\sigma^{2}(z)=-z$, then
$a=\mathrm{id}_{S^{1}}\times\sigma^{2}$, $s=\sigma^{1}\times\mathrm{id}%
_{S^{2}}$. Finally, $r$ can be lifted as follows:
\begin{multline*}
\Psi_{n} \left(  \left(  \theta,\lambda\exp(i\theta) \right)  , \left(
[z,1],\mu\lbrack z^{n},1] \right)  \right)  =\\
\left(  \left(  \theta,\lambda\exp(i\theta) \right)  , \left(  [\exp
(i\theta)z,1],\mu\lbrack\exp(in\theta)z^{n},1] \right)  \right)  ,
\end{multline*}
for all $0\leq\theta\leq2\pi$, $\lambda\in\mathbb{R}$, $\mu\in\mathbb{C}$.
\end{example}

\bigskip

\noindent\textbf{Authors' addresses}

\smallskip

\noindent(J.M.M) \textsc{Institute of Information Security, CSIC, C/ Serrano
144, 28006-Madrid, Spain}

\noindent\textit{E-mail:\/} \texttt{jaime@iec.csic.es}

\medskip

\noindent(M.E.R.M) \textsc{Departament of Applied Mathematics, ETSAM, UPM,
Avda.\ Juan de Herrera 4, 28040-Madrid, Spain}

\noindent\emph{E-mail:\/} \texttt{eugenia.rosado@upm.es}

\medskip

\noindent(I.S.R) \textsc{Departament of Geometry and Topology, Faculty of
Science, University of Granada, Avda.\ Fuentenueva s/n, 18071-Granada, Spain}

\noindent\emph{E-mail:\/} \texttt{ignacios@ugr.es}


\begin{thebibliography}{99}                                                                                               %


\bibitem {ArkMaru}M. Arkowitz, Ken-ichi Maruyama, \emph{Self-homotopy
equivalences which induce the identity on homology, cohomology or homotopy
groups}, Topology Appl.\ \textbf{87} (1998), no.\ 2, 133--154.

\bibitem {Baues}H. J. Baues, \emph{On the group of homotopy equivalences of a
manifold}, Trans.\ Amer.\ Math.\ Soc.\ \textbf{348} (1996), no.\ 12, 4737--4773.

\bibitem {B}D. Betounes, \emph{The geometry of gauge-particle field
interaction: a generalization of Utiyama's theorem}, J.
Geom.\ Phys.\ \textbf{6} (1989), 107--125.

\bibitem {Bl}D. Bleecker, \emph{Gauge Theory and Variational Principles},
Addison-Wesley Publishing Company, Inc., Massachusetts, 1981.

\bibitem {DG}F. Ding, H. Geiges, \emph{The diffeotopy group of} $S^{1}\times
S^{2}$ \emph{via contact topology}, Compositio Math.\ \textbf{146} (2010), 1096--1112.

\bibitem {DL}A. Douady and M. Lazard, \emph{Espaces fibr\'es en alg\`ebres de
Lie et en groupes}, Invent.\ Math.\ \textbf{1} (1966), 133--151.

\bibitem {ET}D.B.A. Epstein and W.P. Thurston, \emph{Transformation Groups and
Natural Bundles}, Proc.\ London Math.\ Soc.\ \textbf{38} (1979) 219--236.

\bibitem {G}C. Godbillon, \emph{\'El\'ements de Topologie Alg\'ebrique},
Hermann, Paris, 1971.

\bibitem {Gr}A. Grothendieck, \emph{A General Theory of Fibre Spaces with
Structure Sheaf}, Report No.\ 4, second printing, University of Kansas, 1965.

\bibitem {GS}V. Guillemin and S.\ Sternberg, \emph{Symplectic techniques in
physics}, Cambridge University Press, Cambridge, U.K., 1984.

\bibitem {H}D. Husemoller, \emph{Fibre Bundles}, Springer Verlag, New York,
Inc., 1994.

\bibitem {KN}S. Kobayashi and K. Nomizu, \emph{Foundations of Differential
Geometry, Volumes I, II}, Interscience Publishers, John Wiley \& Sons, Inc.,
New York, London, 1963, 1969.

\bibitem {KMS}I. Kol\'{a}\v{r}, P.W.\ Michor, J.\ Slov\'ak, \emph{Natural
Operations in Differential Geometry}, Springer-Verlag, Berlin Heidelberg, 1993.

\bibitem {K}J. L. Koszul, \emph{Lectures on Fibre Bundles and Differential
Geometry}, Lectures on Mathematics and Physics no.\ 20, Tata Institute of
Fundamental Research, Bombay, India, 1965.

\bibitem {KM}A. Kriegl, P. W. Michor, \emph{The Convenient Setting of Global
Analysis}, Mathematical Surveys and Monographs, 53. American Mathematical
Society, Providence, RI, 1997.

\bibitem {KMR}A. Kriegl, P. W. Michor, A. Rainer, \emph{An exotic zoo of
diffeomorphism groups on $\mathbb{R}^{n}$}, Ann.\ Global
Anal.\ Geom.\ \textbf{47} (2), 179--222 (2015).

\bibitem {L1}P. Lecomte, \emph{Sur l'alg\`ebre de Lie des sections d'un
fibr\'e en alg\`ebres de Lie}, Ann.\ Inst.\ Fourier (Grenoble) \textbf{30}
(1980), 35--50.

\bibitem {L2}---, \emph{Note on the linear endomorphisms of a vector bundle},
Manuscripta Math.\ \textbf{32} (1980), 231--238.

\bibitem {L3}---, \emph{On the infinitesimal automorphisms of a vector
bundle}, J. Math.\ Pures Appl.\ \textbf{60} (1981), 229--239.

\bibitem {LR}P. B. A. Lecomte, C. Roger, \emph{D\'eformations de l'alg\`ebre
des automorphismes d'un fibr\'e principal},
Proc.\ Kon.\ Nederl.\ Akad.\ Wetensch., Series A, \textbf{92} (1989), 457--463.

\bibitem {LY}A. J. Ledger and K. Yano, \emph{Almost Complex Structures on
Tensor Bundles}, J. Differential Geom.\ \textbf{1} (1967), 355--368.

\bibitem {MM}K. B. Marathe and G. Martucci, \emph{The mathematical foundations
of gauge theories}, North-Holland, Amsterdam, the Netherlands, 1992.

\bibitem {Mi}J. Mickelsson, \emph{Current Algebras and Groups}, Plenum Press,
New York, 1989.

\bibitem {MS}J. W. Milnor and J. D. Stasheff, \emph{Characteristic Classes},
Annals of Mathematics Studies \textbf{76}, Princeton University Press,
Princeton, NJ, 1974.

\bibitem {Mo}K.-P.\ Mok, \emph{Lifts of Vector Fields to Tensor Bundles},
Geom.\ Dedicata \textbf{8} (1979), 61--67.

\bibitem {Moore}N. Moore, \emph{Algebraic Vector Bundles over the }%
$2$\emph{-sphere}, Inventiones Math.\ \textbf{14} (1971), 167--172.

\bibitem {PT}R. S. Palais and Ch.-L.\ Terng, \emph{Natural Bundles have Finite
Order}, Topology \textbf{16} (1977), 271--277.

\bibitem {Pav}P. Pave\v{s}i\'{c}, \emph{On self-maps which induce identity
automorphisms of homology groups}, Glasgow Math.\ J. \textbf{43} (2001),
no.\ 2, 177--184.

\bibitem {Ruber}D. Ruberman, \emph{An Obstruction to Smooth Isotopy in
Dimenson $4$}, Mathematical Research Letters 5 (1998), 743--758.

\bibitem {SS}P. Sankaran, S. Sarkar, \emph{Degrees of maps between Grassmann
manifolds}, Osaka J. Math. \textbf{46} (2009), no.\ 4, 1143--1161.

\bibitem {Sp}E. H. Spanier, \emph{Algebraic Topology}, McGraw-Hill, Inc., 1966.

\bibitem {T}Ch.-L. Terng, \emph{Natural Vector Bundles and Natural
Differential Operators}, Amer.\ J. Math.\ \textbf{100} (1978), 775--828.

\bibitem {YL}K. Yano and A.J. Ledger, \emph{Linear Connections on Tangent
Bundles}, J. London Math.\ Soc.\ \textbf{39} (1964), 495--500.
\end{thebibliography}
\end{document}